\documentclass[10pt]{article}

\usepackage[draft]{graphics}
\usepackage{graphicx}          	 % for eps figures
\usepackage{bm}                 % for bold math symbols
\usepackage{amsmath}             % for \text and such\gstar
\usepackage{amssymb}
\usepackage{amsfonts}             
\usepackage{verbatim}           % useful for commenting out stuff
\usepackage{amsthm}             % does theorems etc. nicely, but

\usepackage{mathtools}
\usepackage{relsize}

\usepackage{makeidx}

\theoremstyle{plain}             % This is the default

\theoremstyle{definition}

\makeatletter
\@addtoreset{equation}{chapter}

\makeatother

\def\protectbold#1{\protect{\boldmath{$#1$}}}

\def\eqref#1{(\ref{#1})}
\def\dsp{\displaystyle}
\def\Frac#1#2{\frac
{
 {\raise.6ex
 \hbox{$\displaystyle#1$}}
}
{
 {\lower.6ex
 \hbox{$\displaystyle#2$}}
 }
}

\numberwithin{equation}{section}

\def\binomial#1#2{
\renewcommand{\arraystretch}{1.0}
\left(
\begin{array}{c} 
\hskip-5pt#1\\
\hskip-5pt#2
\end{array}
\hskip-5pt\right)}

\def\bigOxe{\sqcup \kern-2.3mm \sqcap}

\makeatother

\def\dsp{\displaystyle}
\def\Frac#1#2{\frac
{
 {\raise.6ex
 \hbox{$\displaystyle#1$}}
}
{
 {\lower.6ex
 \hbox{$\displaystyle#2$}}
 }
}

\input epsf

\def\CHF#1#2#3{
{}_1F_1\left(
\begin{array}{c}
\begin{array}{cc} \hskip-10pt#1 \end{array}\\
\begin{array}{c}  \hskip-10pt#2 \end{array}
\end{array}
\hskip-8pt;\,#3
\right)}

\def\CHFs#1#2#3{
{}_1F_1\left({a};{c};{z}\right)
}

\def\intp{\int_0^\infty}

\def\bigO{{\cal O}}
\def\calC{{{\cal C}}}

\def\calL{{{\cal L}}}

\def\ph{{\rm ph}}

\def\tfrac#1#2{{{\lower.6ex
\hbox{$\scriptstyle#1$}}\over 
{\raise.7ex
\hbox{$\scriptstyle#2$}}}}

\def\intp{\int_0^\infty}

\def\calL{{\cal L}}

\def\phase{{\rm ph}}

\def\tfrac#1#2{{{\lower.6ex
\hbox{$\scriptstyle#1$}}\over 
{\raise.7ex
\hbox{$\scriptstyle#2$}}}}

\def\sn{\sum_{n=0}^\infty\,}
\def\sk{\sum_{k=0}^\infty\,}

\def\insil#1{}

\begin{document}
 \title{
Remarks on Slater's asymptotic expansions of Kummer functions for large values of the $a-$parameter
}

\author{
    Nico M. Temme\footnote{Emeritus researcher at Centrum Wiskunde \& Informatica (CWI), 
        Science Park 123, 1098 XG Amsterdam,  The Netherlands}
       \\
        IAA, Abcoude 1391 VD 18,   The Netherlands\\
    {\ }\\
     { \small e-mail: {\tt Nico.Temme@cwi.nl}}\\{\ }\\
     To the memory of Panayiotis D. Siafarikas \\
The man who loved special functions.
}

\date{\ }

\maketitle
\begin{abstract}
In Slater's 1960 standard work on confluent hypergeometric functions, also called Kummer functions, a number of asymptotic expansions  of  these functions can be found.  We summarize  expansions derived from a differential equation  for large values of the $a-$parameter. We show how similar expansions can be derived by using integral representations, and we observe discrepancies with Slater's expansions.
\end{abstract}

\vskip 0.8cm \noindent
{\small
\par\noindent
Keywords \& Phrases: Asymptotic analysis; Kummer functions; confluent hypergeometric functions; Bessel functions.
}

\def\cprime{$'$}

\section{Introduction}\label{sec:intro}

Large parameter problems can be presented in the form of integrals or differential equations, or both, but we also encounter finite sums,  infinite series, difference equations, and implicit algebraic equations. In this paper we  use integral representations of the confluent hypergeometric functions, also called Kummer functions,   and we derive expansions  of  the Kummer functions ${}_1F_1(a;b;z)$ and $U(a,b,z)$ for large positive and negative values of $a$. The expansions are in terms of the modified Bessel functions $I_\nu(z)$ and $K_\nu(z)$, and they are valid for bounded values of $z$ and $b$. 

In the next section we  summarize similar results given in Slater's  standard work \cite{Slater:1960:CHF} on Kummer functions, which results are derived by using Kummer's differential equation. After we have derived our results for ${}_1F_1(a;c;z)$ and $U(a,c,z)$ in subsequent sections by using integral representations,  we can  compare the results of both approaches.  

We observe that Slater's large $a-$expansions  of the $U-$function is  not in agreement with our result for this function, and the question arises which result is correct, and why certain steps leading to wrong results can be explained.

We also mention other large $a-$expansions of Kummer functions available in the literature. For information on the Kummer functions we refer to Adri Olde Daalhuis' chapter on these functions \cite{Olde:2010:CHF}\footnote{See also {http://dlmf.nist.gov/13}.} in the {\em N{IST} {H}andbook of {M}athematical {F}unctions} \cite{Olver:2010:HMF}, and we quote some of the formulas that are relevant in our analysis.

\section{Slater's results}\label{sec:slater}
Slater's expansions for large $a$ are given in \cite[\S4.6.1]{Slater:1960:CHF}, and are in terms of the large parameter 
$a$ written in the form
\begin{equation}\label{eq:sla01}
a=\tfrac14u^2+\tfrac12b,
\end{equation}
where $u>0$ if $a$ and $b$ are real with $a>\frac12b$. Then,\footnote{In an earlier version of this paper the argument of the ${}_1F_1$ was not correct. With thanks to Martin Ehler and Karlheinz Gr{\"o}chenig, who noticed this in 	arXiv:2208.01122.}
\begin{equation}\label{eq:sla02}
\begin{array}{ll}
{e^{-\frac12z^2}z^b\CHF{a}{b}{z^2}=\Gamma(b)u^{1-b}2^{b-1}}\ \times \\[8pt]
\quad\quad
\dsp{\left(zI_{b-1}(uz)\left(\sum_{s=0}^{N-1}\frac{A_s(z)}{u^{2s}}+\bigO\left(\frac{1}{u^{2N}}\right)\right)+\right. }\\[8pt]
\quad\quad\dsp{\left.
\frac{z}{u}I_{b}(uz)\left(\sum_{s=0}^{N-1}\frac{B_s(z)}{u^{2s}}+\frac{z}{1+\vert z\vert}\bigO\left(\frac{1}{u^{2N}}\right)\right)\right)}
\end{array}
\end{equation}
and the incorrect result, see \eqref{eq:kumrem02},
\begin{equation}\label{eq:sla03}
\begin{array}{ll}
\dsp{e^{-\frac12z^2}z^bU\left(a,b,z^2\right)=\frac{2^{2-b}u^{b-1}}{\Gamma(a)}}\ \times \\[8pt]
\quad\quad
\dsp{\left(zK_{b-1}(uz)\left(\sum_{s=0}^{N-1}\frac{A_s(z)}{u^{2s}}+\bigO\left(\frac{1}{u^{2N}}\right)\right)-\right. }\\[8pt]
\quad\quad\dsp{\left.
\frac{z}{u}K_{b}(uz)\left(\sum_{s=0}^{N-1}\frac{B_s(z)}{u^{2s}}+\frac{z}{1+\vert z\vert}\bigO\left(\frac{1}{u^{2N}}\right)\right)\right)}.
\end{array}
\end{equation}
where $I_\nu(z)$ and $K_\nu(z)$ are the modified Bessel functions and the coefficients are given by $A_0=1$ and
\begin{equation}\label{eq:sla04}
\begin{array}{ll}
\dsp{B_s(z)=-\tfrac12A_s^\prime(z)+\int_0^z\left(\tfrac12t^2A_s(t)-\frac{b-\frac12}{t}A_s^\prime(t)\right)\,dt,} \\[8pt]
\dsp{A_{s+1}(z)=\frac{b-\frac12}{z}B_s-\tfrac12B_s^\prime(z)+\int \tfrac12t^2B_s(t)\,dt+K_s},
\end{array}
\end{equation}
and $K_s$ is chosen so that $A_{s+1}(z)\to0$ as $z\to0$. 

In fact\footnote{In Slater's formula (4.6.46) the $x$ should be a $z$.},
\begin{equation}\label{eq:sla05}
\begin{array}{ll}
\dsp{A_0(z)=1,}\\[8pt]
\dsp{B_0(z)=\tfrac16z^3,}\\[8pt]
\dsp{A_1(z)=\tfrac16(b-2)z^2+\tfrac{1}{72}z^6,}\\[8pt]
\dsp{B_1(z)=-\tfrac13b(b-2)z-\tfrac{1}{15}z^5+\tfrac{1}{216}z^{9},}\\[8pt]
\dsp{A_2(z)=-\tfrac{1}{120}(5b-12)(b+2)z^4+\tfrac{1}{6480}(5b-52)z^8+\tfrac{1}{31104}z^{12},}\\[8pt]
\dsp{B_2(z)=\tfrac{1}{90}(5b-12)(b+2)(b+1)z^3-\tfrac{1}{45360}(175b^2-350b-1896)z^7+}\\[8pt]
\quad\quad\quad\quad \dsp{-\tfrac{7}{12960}z^{11}+\tfrac{1}{933120}z^{15}.}
\end{array}
\end{equation}

Slater claims that these expansions are valid uniformly with respect to $z$ in bounded domains. In the next sections we derive expansions of ${}_1F_1(a;b;z)$ and $U(a,b,z)$ for large $a$ and  compare
 these results with Slater's expansions. 
 
 \medskip
 
 {\em It will appear that the expansion for the $U-$function is not correct; see \eqref{eq:kumrem02}.}

 \section{Expansions  for \protectbold{a\to+\infty}}\label{sec:kumlargea}

We derive the expansions of ${}_1F_1(a;b;z)$ and $U(a,b,z)$ by using integral representations. Because we want to compare our results with those of Slater we use for $a$ the form given in \eqref{eq:sla01} and replace $z$ with $z^2$.

\subsection{Expansion of \protectbold{U(a,b,z)}}\label{sec:kumUlargea}
We summarize results from  \cite{Temme:1981:OTE}, but we use the notation used by Slater as in \S\ref{sec:slater}. We start with 
\begin{equation}\label{eq:kumUap01}
U(a,b,z)=\frac{1}{\Gamma(a)}\intp e^{-zt} t^{a-1}(1+t)^{b-a-1}\,dt,
\end{equation}
valid for $\Re a>0$ and $ \Re z>0$. By writing $t/(1+t)=e^{-s}$ we obtain  after a few steps
\begin{equation}\label{eq:kumUap02}
U\left(a,b,z^2\right)=\frac{e^{\frac12z^2}}{\Gamma(a)}\intp e^{-\frac14u^2s-z^2/s} s^{-b}f(s)\,ds,
\end{equation}
where
\begin{equation}\label{eq:kumUap03}
f(s)=e^{z^2\mu(s)}\left(\frac{s/2}{\sinh(s/2)}\right)^b,\quad \mu(s)=\frac{1}{s}-\frac{1}{e^s-1}-\frac12.
\end{equation}
The function $f$ is analytic in the strip $\vert\Im s\vert<2\pi$ and it can be expanded for  $\vert s\vert<2\pi$ into a Maclaurin expansion. We write an expansion with a remainder in the form
\begin{equation}\label{eq:kumUap04}
f(s)=\sum_{k=0}^{K-1} c_k s^k+s^Kr_K(s),\quad K=0,1,2,,\ldots.
\end{equation}
The coefficients $c_k$  are combinations of Bernoulli numbers and Bernoulli polynomials\footnote{http://dlmf.nist.gov/24}. We have
\begin{equation}\label{eq:kumUap05}
\left(\frac{s/2}{\sinh(s/2)}\right)^{b}=e^{\frac12bs}\left(\frac{s}{e^s-1}\right)^b=\sk \frac{B_k^{b}(b/2)}{k!}s^k,
\end{equation}
and 
\begin{equation}\label{eq:kumUap06}
\mu(s)=-\sum_{k=1}^\infty \frac{B_{2k}}{(2k)!}s^{2k-1}.
\end{equation}
The first $c_k$  are 
\begin{equation}\label{eq:kumUap07}
\begin{array}{@{}r@{\;}c@{\;}l@{}}
c_0&=&1,\quad c_1=-\frac{1}{12}z^2,\\[8pt]
c_2&=&\frac{1}{288}\left(z^4-12b\right),\\[8pt]
c_3&=&\frac{z^2}{51840}\left(72+180b-5z^4\right),\\[8pt]
c_4&=&\frac{1}{2488320}\left(5z^8-(288+360b)z^4+864b+2160b^2\right).
\end{array}
\end{equation}

We substitute the expansion in \eqref{eq:kumUap04} into \eqref{eq:kumUap02} and obtain
\begin{equation}\label{eq:kumUap08}
U\left(a,b,z^2\right)=\frac{e^{\frac12z^2}}{\Gamma(a)}\sum_{k=0}^{K-1} c_k\Phi_k+R_K(a,b,z),
\end{equation}
where
\begin{equation}\label{eq:kumUap09}
R_K(a,b,z)=\frac{e^{\frac12z^2}}{\Gamma(a)}\intp e^{-\frac14u^2s-z^2/s} s^{K-b}r_K(s)\,ds,
\end{equation}
and, in terms of the modified Bessel function $K_\nu(z)$,
\begin{equation}\label{eq:kumUap10}
\Phi_k=\intp e^{-\frac14u^2s-z^2/s} s^{k-b}\,ds=2\left(\frac{2z}{u}\right)^{k-b+1}K_{k-b+1}(uz).
\end{equation}
This representation follows from\footnote{http://dlmf.nist.gov/10.32.E10}
\begin{equation}\label{eq:kumUap11}
K_{{\nu}}(z)=\tfrac{1}{2}(\tfrac{1}{2}z)^{\nu}\int _{0}^{\infty}e^{-t-z^2/(4t)}\frac{dt}{t^{{\nu+1}}},\quad \vert\phase\,z\vert<\tfrac14\pi,
\end{equation}
which function is an even function of $\nu$.

In \cite{Temme:1981:OTE} we have constructed a bound for the remainder $R_K$ and we have shown that the sequence $\{\Phi_k\}$ constitutes an asymptotic sequence  for $u\to+\infty$ in the sense that 
\begin{equation}\label{eq:kumUap12}
\frac{\Phi_k}{\Phi_{k-1}}=\bigO\left(\frac{1+uz}{u^2}\right),\quad u\to+\infty,
\end{equation}
uniformly in bounded $b-$intervals and bounded $z-$intervals ($z>0$), but these intervals can be extended to complex domains. This shows the asymptotic nature of the expansion in \eqref{eq:kumUap08}.

We can obtain an expansion with only two Bessel functions by using the recursion
\begin{equation}\label{eq:kumUap13}
K_{\nu+1}(z)=K_{\nu-1}(z)+\frac{2\nu}{z}K_\nu(z).
\end{equation}
and rearranging the expansion. A more direct way follows from writing
\begin{equation}\label{eq:kumUap14}
f(s)=\alpha_0+\beta_0s+s^2g(s),\quad \alpha_0=c_0,\quad \beta_0=c_1.
\end{equation}
Substituting this in \eqref{eq:kumUap02} we obtain after integrating by parts
\begin{equation}\label{eq:kumUap15}
U\left(a,b,z^2\right)=\frac{e^{\frac12z^2}}{\Gamma(a)}\left(
\alpha_0\Phi_0+\beta_0\Phi_1+\frac{1}{u^2}\intp e^{-\frac14u^2s-z^2/s} s^{-b}f_1(s)\,ds\right),
\end{equation}
where
\begin{equation}\label{eq:kumUap16}
f_1(s)=4s^be^{z^2/s}\frac{d}{ds}\left(e^{-z^2/s}s^{2-b}g(s)\right).
\end{equation}

Considering the behavior of $f$ (defined in \eqref{eq:kumUap03}) at infinity, and that of $g_1$ and successive $g_n, f_n$, we observe that $f(s)=\bigO(\exp(-bs/2))$ as $s\to\infty$ when $\Re b<0$; when $\Re b\ge0$, $f$ is bounded. It follows that the integrated term at infinity will vanish if $\Re(u^2s)>0$ and $u$ is large enough. When $u$ is complex, we may turn the path of integration into the complex plane over  an angle $\theta$ with $\vert\theta\vert<\frac12\pi$. This  is possible if $-\pi+\delta\le  \ph(u^2)\le\pi-\delta$, with $\delta$ a small positive number. 

At the origin $f$ and $g_1$ (and successive $f_n$ and $g_n$) are analytic, and the integrated term will vanish if $\Re(z^2/s)>0$. Again, when $z$ is complex, we may achieve this by integrating from the origin in a suitable direction, and deform the contour to get a suitable direction at infinity.

The integration by parts procedure can be continued, and we obtain
\begin{equation}\label{eq:kumUap17}
\begin{array}{ll}
\dsp{U\left(a,b,z^2\right)=\frac{e^{\frac12z^2}}{\Gamma(a)}\left(\Phi_0\sum_{n=0}^{N-1} \frac{\alpha_n}{u^{2n}}+\Phi_1\sum_{n=0}^{N-1} \frac{\beta_n}{u^{2n}}\right. +}\\[8pt]
\quad\quad\quad\quad
\dsp{\left.\frac{1}{u^{2N}}\intp e^{-\frac14u^2s-z^2/s} s^{-b}f_N(s)\,ds\right)},
\end{array}
\end{equation}
where the $\Phi_k$ are defined in \eqref{eq:kumUap10} and $\alpha_n, \beta_n, f_n$ follow from the recursive scheme
\begin{equation}\label{eq:kumUap18}
\begin{array}{@{}r@{\;}c@{\;}l@{}}
f_n(s)&=&\alpha_n+\beta_ns+s^2g_n(s), \\[8pt]
f_{n+1}(s)&=&\dsp{4 s^be^{z^2/s}\frac{d}{ds}\left(e^{-z^2/s}s^{2-b}g_n(s)\right)},
\end{array}
\end{equation}
with $f_0=f$.

We can express the coefficients $\alpha_n$ and $\beta_n$ in terms of the $c_k$ used in \eqref{eq:kumUap04}. We write
\begin{equation}\label{eq:kumUap19}
f_n(s)=\sk c_k^{(n)} s^k, \quad c_k^{(0)}=c_k,
\end{equation}
and after substituting this into \eqref{eq:kumUap18} we find for the coefficients $c_k^{(n)}$ the recursion
\begin{equation}\label{eq:kumUap20}
c_0^{(n+1)} =4z^2c_2^{(n)} ,\quad
c_k^{(n+1)} =4\left(z^2c_{k+2}^{(n)}+(1-b+k)c_{k+1}^{(n)}\right),
\end{equation}
where $k\ge1$ and $n\ge0$. The first coefficients are
\begin{equation}\label{eq:kumUap21}
\begin{array}{@{}r@{\;}c@{\;}l@{}}
\alpha_0&=&1,\quad  \beta_0=c_1,\\[8pt]
\alpha_1&=&4z^2c_2,\quad  \beta_1=4z^2c_3+4(2-b)c_2\\[8pt]
\alpha_2&=&4z^2(4z^2c_4+4(3-b)c_3),\\[8pt]
\beta_2&=&16z^4c_5+32z^2((3-b)c_4+16(b-2)(b-3)c_3.
\end{array}
\end{equation}
In general, for $\alpha_n$ we need $c_{n+1},\cdots,c_{2n}$ and for $\beta_n$ we need $c_{n+1},\cdots,c_{2n+1}$.

To compare the expansion in \eqref{eq:kumUap17} with Slater's expansion in \eqref{eq:sla03}, we observe first that $\Phi_1=2(2z/u)^{2-b}K_{2-b}(uz)$, and we use the relation in \eqref{eq:kumUap13} to rearrange our expansion. This gives (we have used $K_\nu(z)=K_{-\nu}(z)$)
\begin{equation}\label{eq:kumUap22}
\begin{array}{ll}
\dsp{e^{-\frac12z^2}z^bU\left(a,b,z^2\right)=\frac{2^{2-b}u^{b-1}}{\Gamma(a)}}\ \times \\[8pt]
\quad\quad
\dsp{\left(zK_{b-1}(uz)\sum_{n=0}^{N-1}\frac{a_n(z)}{u^{2n}}-\frac{z}{u}K_{b}(uz)\sum_{n=0}^{N-1}\frac{b_n(z)}{u^{2n}}\right. +}\\[8pt]
\quad\quad\quad\quad\quad\quad
\dsp{\left.\frac{2^{b-2}z^bu^{1-b}}{u^{2N}}\intp e^{-\frac14u^2s-z^2/s} s^{-b}f_N(s)\,ds\right)},
\end{array}
\end{equation}
where 
\begin{equation}\label{eq:kumUap23}
\begin{array}{@{}r@{\;}c@{\;}l@{}}
a_0(z)&=&1,\quad  a_n(z)= \alpha_n+4(1-b)\beta_{n-1},\quad n\ge1,\\[8pt]
b_n(z)&=&-2z\beta_n,\quad n\ge0.
\end{array}
\end{equation}
This gives the first coefficients
\begin{equation}\label{eq:kumUap24}
\begin{array}{ll}
\dsp{a_0(z)=1,}\\[8pt]
\dsp{b_0(z)=\tfrac16z^3,}\\[8pt]
\dsp{a_1(z)=\tfrac16(b-2)z^2+\tfrac{1}{72}z^6,}\\[8pt]
\dsp{b_1(z)=-\tfrac13b(b-2)z-\tfrac{1}{15}z^5+\tfrac{1}{1296}z^{9},}\\[8pt]
\dsp{a_2(z)= -\tfrac23b(b-1)(b-2)-\tfrac{1}{120}(b+2)(5b-12)z^4+}\\[8pt]
\quad\quad\quad\quad \dsp{\tfrac{1}{6480}(5b-52)z^8+\tfrac{1}{31104}z^{12},}\\[8pt]
\dsp{b_2(z)=
-\tfrac{1}{90}(5b+2)(b-3)(b-4)z^3-\tfrac{1}{45360}(175b^2-350b-1896)z^7+}\\[8pt]
\quad\quad\quad\quad \dsp{-\tfrac{7}{12960}z^{11}+\tfrac{1}{933120}z^{15}.}
\end{array}
\end{equation}

When we compare these coefficients with the ones in Slater's expansion of the $U-$function given in \eqref{eq:sla05} we see differences in $a_2(z)$ and $b_2(z)$. In particular, the condition $A_n(0)=0$  ($n\ge1$) used in the construction of Slater's coefficients is not showing in our $a_2(z)$.

\subsection{Expansion of \protectbold{{}_1F_1(a;b;z)}}\label{sec:kumFlargea}
For an expansion of the $F-$function we start with the integral\footnote{http://dlmf.nist.gov/13.4.ii}
\begin{equation}\label{eq:kumFap01}
\CHF{a}{b}{z}=\frac{\mathop{\Gamma(b)\Gamma\/}\nolimits\!\left(1+a-b\right)}{2\pi i\mathop{\Gamma\/}\nolimits\!\left(a\right)}\int _{0}^{{(1+)}}e^{{zt}}t^{{a-1}}{(t-1)^{{b-a-1}}}\,dt,
\end{equation}
where $\Re a>0$ and  $b-a\ne 1,2,3,\ldots$.
The contour can be the circle $\vert t-1\vert=1$. The transformation $t=s/(s-1)$ transforms this circle into itself. To verify this we write $s=t/(t-1)$. With $t=1+e^{i\theta}$, $\theta\in[0,2\pi)$, we obtain $s=1+e^{-i\theta}$. The result of the substitution is
\begin{equation}\label{eq:kumFap02}
\CHF{a}{b}{z}=\frac{\mathop{\Gamma(b)\Gamma\/}\nolimits\!\left(1+a-b\right)}{2\pi i\mathop{\Gamma\/}\nolimits\!\left(a\right)}\int _\calC e^{{zs/(s-1)}}s^{{a-1}}{(s-1)^{{b}}}\,ds,
\end{equation}
where $\calC$ is the circle $\vert s-1\vert=1$.

Next we take $s=e^w$. With $s=1+e^{i\theta}$, $\theta\in[0,2\pi)$, we see that the circle $\calC$ is described by 
\begin{equation}\label{eq:kumFap03}
w=\sigma+i\tau,\quad \sigma=\ln(2\cos \tau),\quad -\tfrac12\pi<\tau<\tfrac12\pi.
\end{equation}
After some manipulations we obtain
\begin{equation}\label{eq:kumFap04}
\CHF{a}{b}{z^2}=\frac{\mathop{\Gamma(b)\Gamma\/}\nolimits\!\left(1+a-b\right)e^{\frac12z^2}}{\mathop{\Gamma\/}\nolimits\!\left(a\right)\,2\pi i}\int _\calL e^{\frac14u^2s + z^2/s}s^{{-b}}f(-s)\,ds,
\end{equation}
where $a=\frac14u^2+\frac12b$ (as in \eqref{eq:sla01}), $f$ is the same as in \eqref{eq:kumUap03} and $\calL$ can be taken as a loop around the negative axis that encircles the origin in a positive (anti-clockwise) direction. Below and above the branch cut along the negative axis the phase of $s$ is $-\pi$ and $+\pi$, respectively. This representation is valid for all complex $z$ and $\Re(a+b)>0$.

Upon substituting the expansion in  \eqref{eq:kumUap04} we obtain
\begin{equation}\label{eq:kumFap05}
\frac{1}{\Gamma(b)}\CHF{a}{b}{z^2}=\frac{\Gamma(1+a-b)e^{\frac12z^2}}{\Gamma(a)}\sum_{k=0}^{K-1}(-1)^k  c_k\Psi_k+S_K(a,b,z),
\end{equation}
where
\begin{equation}\label{eq:kumFap06}
S_K(a,b,z)=(-1)^K\frac{\Gamma(1+a-b) e^{\frac12z}}{\Gamma(a)\,2\pi i}\int_\calL e^{\frac14u^2s+z^2/s} s^{K-b}r_K(-s)\,ds,
\end{equation}
and, in terms of the modified Bessel function $I_\nu(z)$,
\begin{equation}\label{eq:kumFap07}
\Psi_k=\frac{1}{2\pi i}\int_\calL e^{\frac14u^2s+z^2/s} s^{k-b}\,ds =\left(\frac{2z}{u}\right)^{k+1-b}I_{b-k-1}(uz).
\end{equation}
This representation follows from\footnote{http://dlmf.nist.gov/10.9.E19}
\begin{equation}\label{eq:kumFap08}
J_{{\nu}}(z)=\frac{(\tfrac{1}{2}z)^{\nu}}{2\pi i}\int _{{-\infty}}^{{(0+)}}e^{t- z^{2}/(4t)}\frac{dt}{t^{{\nu+1}}},
\end{equation}
with $z$ replaced with $e^{\frac12\pi i}z$.

In the above results we can give $z$ any finite complex value, and we require $\Re a>0$, $1+a-b\ne0,-1,-2,\ldots$.
For $b=0,-1,-2,\ldots$, the left-hand side of \eqref{eq:kumFap05} can be interpreted by using 
\begin{equation}\label{eq:kumFap09}
\lim_{b\to-m}\frac{1}{\Gamma(b)}\CHF{a}{b}{z}=\frac{(a)_{m+1}\,z^{m+1}}{(m+1)!}\CHF{a+m+1}{m+2}{z}.
\end{equation}

The expansion in \eqref{eq:kumFap05} can be written in the form with two Bessel functions. We need the relation
\begin{equation}\label{eq:kumFap10}
I_{\nu-1}(z)=I_{\nu+1}(z)+\frac{2\nu}{z}I_\nu(z), 
\end{equation}
and an integration by parts procedure as used for the $U-$function gives a form comparable with Slater's result in \eqref{eq:sla02}. In this way we obtain the result written in the form of \eqref{eq:sla02}
\begin{equation}\label{eq:kumFap11}
\begin{array}{ll}
\dsp{e^{-\frac12z^2}z^b\CHF{a}{b}{z^2}=\frac{\Gamma(b)\Gamma(1+a-b)}{\Gamma(a)}u^{b-1}2^{1-b}}\ \times \\[8pt]
\quad\quad
\dsp{\left(zI_{b-1}(uz)\sum_{n=0}^{N-1}\frac{a_n(z)}{u^{2n}}+\frac{z}{u}I_{b}(uz)\sum_{n=0}^{N-1}\frac{b_n(z)}{u^{2n}}\right. +}\\[8pt]
\quad\quad
\dsp{\left.\frac{2^{b-1}z^bu^{1-b}}{u^{2N}\,2\pi i}\int_{\calL} e^{\frac14u^2s+z^2/s} s^{-b}f_N(-s)\,ds\right)},
\end{array}
\end{equation}
where the coefficients $a_n(z), b_n(z)$ and the functions $f_n$ are the same as for the expansion of the $U-$function in 
\eqref{eq:kumUap22}.

We cannot yet compare  this result with Slater's result in \eqref{eq:sla02}, because of the ratio of the gamma functions with large parameter in our results. We should expand this ratio and multiply this expansion with the ones in \eqref{eq:kumFap11}. We have
\begin{equation}\label{eq:kumFap12}
\frac{\Gamma(1+a-b)}{\Gamma(a)}=\frac{\Gamma\left(1+\frac14u^2-\frac12b\right)}{\Gamma\left(\frac14u^2+\frac12b\right)}\sim \left(\frac{u}{2}\right)^{2-2b}\sn \frac{d_n}{u^{2n}},
\end{equation}
as $a\to\infty$. All  coefficients $d_{2n+1}$ vanish and the first even indexed coefficients are
\begin{equation}\label{eq:kumFap13}
\begin{array}{@{}r@{\;}c@{\;}l@{}}
d_0&=&1,\quad d_2=\frac{2}{3}(b-2)(b-1)_2,\\[8pt]
d_4&=&\frac{2}{45}\left(5b^2-22b+24\right)(b-1)_4,\\[8pt]
d_6&=&\frac{4}{2835}\left(35b^3-252b^2+604b-480\right)(b-1)_6,\\[8pt]
d_8&=&\frac{2}{42525}\left(175b^4-1820b^3+7124b^2-12400b+8064\right)(b-1)_8.
\end{array}
\end{equation}
See \cite[\S3.6.2]{Temme:1996:SFA} for expansions of this type.

When we perform the multiplications of the series we obtain Slater's expansion given in \eqref{eq:sla02}.

\subsection{Remarks on both methods}\label{sec:kumrems}
Slater's expansions in \S\ref{sec:slater} are based on Olver's method for differential equations; see \cite[Chapter 12]{Olver:1974:ASF} (Slater has referred to earlier papers by Olver). This method is very powerful, it gives expansions valid in large domains of the parameters and recurrence relations for the coefficients. Also, the method provides realistic error bounds for remainders in the expansions. In the case of the Kummer functions, the expansions are first given for two linear independent solutions of Kummer's differential equation
\begin{equation}\label{eq:kumrem01}
zw^{\prime\prime}+(b-z)w^\prime-aw=0.
\end{equation}
Then the expansions of ${}_1F_1(a;b;z)$ and $U(a,b,z)$ follow from linear combinations of these solutions, and the coefficients in these combinations follow from certain known limiting forms of the Kummer functions (in the present case for $z\to0$). 

On the other hand, when the recurrence relations for the coefficients in the expansions are derived, these recursions usually include constants of integration. In the present case these are the quantities $K_s$ used by Slater in \eqref{eq:sla04}. 
A certain choice of these constants generates a formal solution of the differential equation.

These two steps have to be taken into account when constructing the expansions of the functions ${}_1F_1(a;b;z)$ and $U(a,b,z)$, and it appears that Slater has not used the correct steps for the $U-$function.

When working with integrals these difficulties are not present: we always start with a representation of the function to be considered. All right, we can usually not construct recurrence relations for the coefficients, and the construction of error bounds or estimates for remainders is more difficult, but there will never be a misunderstanding about the correct form of the expansions. 

We can repair Slater's expansion by dividing the series by the series in \eqref{eq:kumFap12}, which gives the expansion given in \eqref{eq:kumUap17}. We can also repair by including a ratio of gamma functions in the representation in \eqref{eq:sla03}, and modify powers of $2$ and $u$. That is, we have

\begin{equation}\label{eq:kumrem02}
\begin{array}{ll}
\dsp{e^{-\frac12z^2}z^bU\left(a,b,z^2\right)=\frac{2^{b}u^{1-b}}{\Gamma(1+a-b)}}\ \times \\[8pt]
\quad\quad
\dsp{\left(zK_{b-1}(uz)\left(\sum_{s=0}^{N-1}\frac{A_s(z)}{u^{2s}}+\bigO\left(\frac{1}{u^{2N}}\right)\right)-\right. }\\[8pt]
\quad\quad\dsp{\left.
\frac{z}{u}K_{b}(uz)\left(\sum_{s=0}^{N-1}\frac{B_s(z)}{u^{2s}}+\frac{z}{1+\vert z\vert}\bigO\left(\frac{1}{u^{2N}}\right)\right)\right)}.
\end{array}
\end{equation}

In the present case the construction of an error bound of the expansion given for the $U-$function given in \eqref{eq:kumUap22} is rather easy when we assume that we have a bound of $f_N(s)$ for $s\ge0$. When $\Re z^2\ge0$, this bound may be independent of $z$, which shows the uniform character of this expansion with respect to $z$. More details on construction of a bound for the remainder in the expansion given in \eqref{eq:kumUap08} can be found in \cite{Temme:1981:OTE}. For the expansions of the $F-$function these bounds should be obtained from complex contours of integration, which is a more difficult matter; again, see \cite{Temme:1981:OTE}.

\subsection{Other forms of the expansions for large \protectbold{a}}\label{sec:kumother}
We have already two forms of the expansions: one with a series of Bessel functions, and one with only two Bessel functions. These forms are valid for bounded and even small values of $z$. When $z$ is such that $uz\to\infty$ we can expand the Bessel functions and use the well-known expansions of these functions for large argument; see \cite[\S10.40]{Olver:2010:BFS}. In this way we can construct  an expansion in terms of elementary functions.

It may also be convenient to have expansions that show the parameter $a$ explicitly as the large parameter, and not the parameter $u$ as in Slater's expansions. Slater's form has some advantages because the coefficients  are simpler than expansions in terms of negative powers of $a$. This is the approach used in \cite{Temme:1981:OTE}, and it gives expansions with a series of Bessel functions. The integration by parts procedure used in the present paper is easy to modify for obtaining expansions with only two Bessel functions. And when we expand the Bessel functions we find expansions with series in negative powers of $a$. For the $F-$function  Perron \cite{Perron:1921:UDV} has given an expansion in terms of elementary functions. For a more recent publication, see \cite{Borwein:2008:ELA}, where an expansion is given for the Laguerre polynomials for large degree. That expansion can also be used for the $F-$function as $a\to-\infty$, because
\begin{equation}\label{eq:kumrem03}
L_n^{(\alpha)}(z)=\binomial{n+\alpha}{n}\CHF{-n}{\alpha+1}{z}.
\end{equation}
Szeg{\H o} \cite[\S8.22, \S8.72, Problemº 46]{Szego:1975:OP} has suggested several methods for these polynomials.

Finally, when $z$ is complex, it may be convenient to consider expansions of the $F-$function in terms of the $J-$Bessel function by using\footnote{http://dlmf.nist.gov/10.27.E6}
\begin{equation}\label{eq:kumrem04}
I_\nu(z)=e^{\mp\frac12\nu\pi i}J_\nu\left(z e^{\pm\frac12\pi i}\right),\quad -\pi\le \pm \phase\,z\le\tfrac12\pi.
\end{equation}

\section{Expansions for \protectbold{a\to-\infty}}\label{sec:kumnega}
For this case we use relations between the Kummer functions and the results for $a\to+\infty$.

\subsection{Expansion of \protectbold{{}_1F_1(a;b;z)}}\label{sec:kumFnega}
This case has not been considered in Slater's book, but we can use the results for $a\to+\infty$ by using  the relation
\begin{equation}\label{eq:kumnega01}
\CHF{a}{b}{z}=e^z\CHF{b-a}{b}{-z}.
\end{equation}
We take this time
\begin{equation}\label{eq:kumnega02}
a=-\tfrac14u^2+\tfrac12b.
\end{equation}
Then $b-a=\frac14u^2+\frac12b$ and we have
\begin{equation}\label{eq:kumnega03}
\CHF{a}{b}{-z^2}=e^{-z}\CHF{\frac14u^2+\frac12b}{b}{z^2},
\end{equation}
For the $F-$function in the right-hand side we can use the results of \S\ref{sec:kumFlargea}. The explicit result for the left-hand side is
\begin{equation}\label{eq:kumnega04}
\begin{array}{ll}
\dsp{e^{\frac12z^2}z^b\CHF{-\frac14u^2+\frac12b}{b}{-z^2}=\frac{\Gamma(b)\Gamma\left(1+\frac14u^2-\frac12b\right)}{\Gamma\left(\frac14u^2+\frac12b\right)}u^{b-1}2^{1-b}}\ \times \\[8pt]
\quad\quad
\dsp{\left(zI_{b-1}(uz)\sum_{n=0}^{N-1}\frac{a_n(z)}{u^{2n}}+\frac{z}{u}I_{b}(uz)\sum_{n=0}^{N-1}\frac{b_n(z)}{u^{2n}}\right. +}\\[8pt]
\quad\quad\quad\quad\quad\quad
\dsp{\left.\frac{2^{b-1}z^bu^{1-b}}{u^{2N}\,2\pi i}\int_{\calL} e^{\frac14u^2s+z^2/s} s^{-b}f_N(-s)\,ds\right)},
\end{array}
\end{equation}
where the coefficients $a_n(z), b_n(z)$ and the functions $f_n$ are the same as for the expansion of the $U-$function in 
\eqref{eq:kumUap22}.

\subsection{Expansion of \protectbold{U(a,b,z)}}\label{sec:kumUnega}
Also in this case we can use connection formulas. We have
\begin{equation}\label{eq:kumnega05}
U(a,b,z)=\frac{\Gamma(1-b)}{\Gamma(a-b+1)}\CHF{a}{b}{z}+
\frac{\Gamma(b-1)}{\Gamma(a)}z^{1-b}\CHF{a-b+1}{2-b}{z},
\end{equation}
when $b$ is not an integer, and
\begin{equation}\label{eq:kumnega06}
\frac1{\Gamma(b)}\CHF{a}{b}{z}=\frac{e^{\mp \pi ia}}{\Gamma(b-a)}U(a,b,z)+
\frac{e^{\pm \pi i(b-a)}}{\Gamma(a)}\,e^zU\left(b-a,b,ze^{\pm
\pi i}\right).
\end{equation}

The first form is useful because we have a real representation,  but, although the $U-$function is well-defined for integer values of $b$, a nasty limiting procedure is needed in that case. The best approach is using the second form with $a=-\frac14u^2 +\frac12b$. Then we can use the expansion given in \eqref{eq:kumnega04} and for 
$U\left(b-a,b,ze^{\pm \pi i}\right)=U\left(\frac14u^2+\frac12b,b,ze^{\pm \pi i}\right)$ the expansion given in \eqref{eq:kumUap22}. Observe that the $K-$Bessel function have arguments $uze^{\pm\frac12\pi i}$, which can be expressed in terms of the ordinary Bessel functions $J_\nu(z)$ and $Y_\nu(z)$; see \cite[\S10.27]{Olver:2010:BFS}.

\section*{Acknowledgements}
The author thanks the referee for useful comments on an earlier version of the paper. He acknowledges support from  {\emph{Ministerio de Ciencia e Innovaci\'on}, Spain}, 
project MTM2009-11686. 

\bibliographystyle{plain}

\def\cprime{$'$} \def\cprime{$'$}

\end{document}